\documentclass[12pt]{article}

\let\Horig\H

\usepackage{amsmath,amscd,amsthm,makeidx,graphicx,url,bm,bbm,mathrsfs,palatino,bm,mathtools,latexsym}
\usepackage[bbgreekl]{mathbbol}

\usepackage{amssymb}
\usepackage{a4wide,xy,graphicx}
\usepackage{tikz-cd}
\usepackage{enumerate,comment}
\usepackage{titlesec}
\usepackage{tikz-cd}

\usepackage[top=2.1cm, bottom=2.1cm, left=2.05cm, right=2.05cm]{geometry}

\numberwithin{equation}{section}

\titleformat{\paragraph}[runin]
{\normalfont\normalsize\bfseries}{\theparagraph}{1em}{}
\titlespacing*{\paragraph}
{0pt}{3.25ex plus 1ex minus .2ex}{1.5ex plus .2ex}

\usepackage[
breaklinks=true,
pdftitle={Commutative avatars of representations of semisimple Lie groups},
pdfauthor={Tam\'as Hausel}]{hyperref}
\usepackage{color}
\definecolor{tocolor}{rgb}{.1,.1,.5}
\definecolor{urlcolor}{rgb}{.2,.2,.6}
\definecolor{linkcolor}{rgb}{.1,.4,.6}
\definecolor{citecolor}{rgb}{.6,.3,.1}
\hypersetup{colorlinks=true, urlcolor=urlcolor, linkcolor=linkcolor, citecolor=citecolor}

\makeindex

\xyoption{all}

\newtheorem{theorem}{Theorem}

\newtheorem{corollary}[theorem]{Corollary}
\newtheorem{conjecture}[theorem]{Conjecture}

\theoremstyle{remark}
\newtheorem{remark}[theorem]{Remark}

\theoremstyle{definition}

\numberwithin{theorem}{section}

\newcounter{margin}
	{\end{itshape}  \bigskip}

\def\beq{\begin{eqnarray}}
	\def\eeq{\end{eqnarray}}
\def\bes{\begin{eqnarray*}}
	\def\ees{\end{eqnarray*}}

\def\s{\mathfrak{s}}

\DeclareMathOperator{\Spec}{Spec} 
\DeclareMathOperator{\Hom}{Hom}

\DeclareMathOperator{\rank}{rank}

\def\C{\mathbb{C}}

\def\M{{\mathcal{M}}}

\def\calQ{{\mathcal{Q}}}
\def\calA{{\mathcal{A}}}

\def\calP{\mathcal{P}}

\def\P{\mathbb{P}}
\def\rP{\rm{P}}
\newcommand{\T}{{\mathrm{T}}}

\def\t{\mathfrak{t}}
\def\calC{{\mathcal C}}
\def\calG{{\mathcal G}}
\def\calB{{\mathcal B}}
\def\calM{{\mathcal M}}
\def\calZ{{\mathcal Z}}

\def\gl{{\mathfrak g\mathfrak l}}
\def\sl{{\mathfrak s\mathfrak l}}
\def\so{{\mathfrak s\mathfrak o}}

\newcommand{\nc}{\newcommand}

\def\calZ{\mathcal{Z}}

\newcommand{\Ad}{\textnormal{Ad}}
\newcommand{\g}{\mathfrak{g}}
\usepackage{leftidx}

\newcommand{\Gr}{\textnormal{Gr}}

\newcommand{\reg}{\textnormal{reg}}

\newcommand{\End}{\textnormal{End}}
\newcommand{\tr}{\textnormal{tr}}

\nc{\op}[1]{\mathop{\mathchoice{\mbox{\rm #1}}{\mbox{\rm #1}}
		{\mbox{\rm \scriptsize #1}}{\mbox{\rm \tiny #1}}}\nolimits}
\nc{\al}{\alpha}

\nc{\ep}{\varepsilon} 
\nc{\ga}{\gamma} 
\nc{\Ga}{\Gamma}
\nc{\la}{\lambda} 
\nc{\La}{\Lambda} 
\nc{\si}{\sigma}
\nc{\Sig}{{\Gamma}} 
\nc{\Om}{\Omega} 
\nc{\om}{\omega}

\nc{\SL}{\mathrm{SL}} 
\nc{\GL}{\mathrm{GL}} 
\nc{\SO}{\mathrm{SO}} 
\nc{\Sp}{\mathrm{Sp}} 
\nc{\PSp}{\mathrm{PSp}}

\nc{\PGL}{\mathrm{PGL}}
\nc{\G}{\mathrm{G}}

\nc{\W}{\mathrm{W}}
\nc{\Lg}{\mathrm{L}}
\nc{\Pg}{\mathrm{P}}
\nc{\calL}{{\mathcal L}}
\nc{\Sym}{{\rm Sym}}

\renewcommand{\H}{\mathrm{H}}

\nc{\Frob}{\mathrm{Frob}}
\nc{\spec}{{\rm Spec}}
\def\PP {\mathbb{P}}

\nc{\rN}{\mathrm{N}}

\usepackage{longtable}

\DeclareMathOperator{\Lie}{Lie}

\nc{\cpt}{{\op{cpt}}} \nc{\Dol}{{\op{Dol}}} \nc{\DR}{{\op{DR}}}
\nc{\B}{{\op{B}}} \nc{\Triv}{\op{Triv}} \nc{\Hod}{{\op{Hod}}}
\nc{\Log}{{\op{Log}}} \nc{\Exp}{{\op{Exp}}} \nc{\Est}{E_{\op{st}}}
\nc{\Hst}{H_{\op{st}}} \nc{\Left}[1]{\hbox{$\left#1\vbox to
		10.5pt{}\right.\nulldelimiterspace=0pt \mathsurround=0pt$}}
\nc{\Right}[1]{\hbox{$\left.\vbox to
		10.5pt{}\right#1\nulldelimiterspace=0pt \mathsurround=0pt$}}
\nc{\LEFT}[1]{\hbox{$\left#1\vbox to
		15.5pt{}\right.\nulldelimiterspace=0pt \mathsurround=0pt$}}
\nc{\RIGHT}[1]{\hbox{$\left.\vbox to
		15.5pt{}\right#1\nulldelimiterspace=0pt \mathsurround=0pt$}}

\nc{\bee}{{\bf E}} 

\title{Commutative avatars of representations of semisimple Lie groups}

\author{ Tam\'as Hausel
\\ {\it IST Austria} 
\\{\tt tamas.hausel@ist.ac.at}   }

\begin{document}

\maketitle

\begin{abstract}
 Here we announce the construction and  properties  of a  big commutative subalgebra of the Kirillov algebra, called big algebra, attached to a finite dimensional irreducible representation of a complex semisimple Lie  group. They are commutative finite flat algebras over the cohomology of the classifying space of the group. They are isomorphic with the equivariant intersection cohomology of affine Schubert varieties, endowing the latter with a new ring structure. Study of the finer aspects of  the structure of the big algebras will also furnish the stalks of the intersection cohomology with ring structure, thus ringifying Lusztig's $q$-weight multiplicity polynomials i.e. certain affine Kazhdan-Lusztig polynomials.  

\end{abstract}

\maketitle

\section{Kirillov and medium algebras}

{L}et $\G$ be a connected complex semisimiple Lie group with Lie algebra $\g$, which we identify with $\g\cong \g^*$ using the Killing form. Let $\mu\in \Lambda^+(\G)$ be a dominant weight, and let $\rho_\mu: \G\to \GL(V^\mu)$ and $\varrho_\mu:=\Lie(\rho_\mu): \g\to \gl(V^\mu)\cong\End(V^\mu)$ be the corresponding complex highest weight representations of the group and its Lie algebra. Using the natural action of $\G$ on the symmetric  algebra $S^*(\g)$  and on the endomorphism algebra $\End(V^\mu)$  Kirillov \cite{kirillov0} introduced  $$\calC^\mu(\g)=\calC^\mu:=(S^*(\g)\otimes \End(V^\mu))^\G\cong 
Maps(\g\cong\g^*\to \End(V^\mu))^\G$$
which we call  {\em (classical) Kirillov algebra}.

Kirillov's motivation for the introduction of $\calC^\mu$ was to understand weight multiplicities of a maximal torus ${\rm T}\subset \G$. For example he proved in \cite[Theorem S]{kirillov0} that $\calC^\mu$ is commutative  if and only if $V^\mu$ is weight multiplicity free. This means that for all $\lambda\in \Lambda=\Hom({\rm T}, \C^\times)$ the weight space $V^\mu_\lambda$ is at most one dimensional. We will see below, that the big commutative subalgebras of the Kirillov algebra we will introduce in this paper will induce in Corollary~\ref{multiplicity} a graded ring structure on multiplicity spaces. 

The Kirillov algebra $\calC^\mu$ is an associative, graded $H_\G^{2*}:=S^*(\g)^\G\cong\C[\g^*]^\G\cong \C[\g]^\G$-algebra. The grading is induced from the usual grading on $S^*(\g)$ and the commutative graded $\C$-algebra $H_\G^{2*}$ acts by scalar multiplication.

 We fix a principal $\sl_2$-subalgebra $\langle e,f,h\rangle\subset \g$ so that we get a section of $
\chi:\g\to \g/\!/\G$, the {\em Kostant section} $\s:=e+\g_f\subset \g^{reg}$, in particular $\s\cong \g/\!/ \G$. Moreover $\s\subset \g^{reg}$ contains only regular elements, i.e. ones with smallest dimensional centralizers, and $\s$ intersects every $\G$-orbit of $\g^{reg}$ in exactly one point. Because the codimension of $\g\setminus \g^{reg}$ in $\g$ is $3$ we can identify \begin{align} \label{kostant}\calC^\mu\cong Maps(\g^{\reg}\to \End (V^\mu))^\G\cong Maps\left(f:\s\to \End(V^\mu) \mid f(x)\in (\End (V^\mu))^{\G_x}\right). \end{align}
We can restrict any subalgebra $\calA\subset \calC^\mu$  to $x\in \g$ to get the finite matrix algebra\begin{align} \label{restrict}\calA_x:=\{f(x) \mid f\in \calA \}\subset (\End(V^\mu))^{\G_x}.
\end{align}

We will denote the one-parameter subgroup $\H^z:\C^\times\to \G_{ad}=\G/Z(\G)$ integrating $\langle h\rangle\subset \g$. Then $\Ad(\H^z)e=z^{-1}e$ and so the $\C^\times$-action \begin{align} \label{action} \begin{array}{ccc}\C^\times \times \g &\to &\g \\ (z,x)&\mapsto& z\cdot x:=\Ad(\H^z)zx \end{array}\end{align} on $\g$ preserves $e$ and $\g_f$ and thus  the Kostant section $\s$, and  induces the grading on $\calC^\mu$ in \eqref{kostant}.

The most important element of $\calC^\mu$, called the {\em small operator} is given by \begin{align} \label{small}\begin{array}{cccc}M_1: &\g&\to& \End(V^\mu)\\ &A &\mapsto &\varrho_\mu(A)\end{array}.\end{align}

More generally we will have an element of the Kirillov algebra from any $\G$-equivariant polynomial map $F:=\g\to\g$ by \begin{align} \label{mediumdef}\begin{array}{cccc}M_F: &\g&\to& \End(V^\mu)\\ &A &\mapsto &\varrho_\mu(F(A))\end{array}.\end{align} For an invariant polynomial $p\in \C[\g]^\G$ we can define 
its derivative $dp:\g\to \g^*\cong \g.$ As $dp$ is automatically $\G$-equivariant we have the operator  $M_{dp}$ from \eqref{mediumdef},  which we call a {\em medium operator} corresponding to $p\in \C[\g]^\G$. For example we have the small operator of \eqref{small} $M_1=M_{d\kappa/2}$ where $\kappa$, the Killing form, is thought of as a degree $2$ invariant polynomial. 
In general we will fix a generating set $\C[\g]^\G\cong \C[p_1,\dots,p_r]$ of homogeneous invariant polynomials $p_i\in \C[\g]^\G$ of degree $d_i$, s.t. $d_1\leq \dots \leq d_r$, where $r=\rank(\G)$. Then we also denote $M_i:=M_{dp_i}$. We will  arrange that $p_1=\kappa/2$ so that $M_1=M_{dp_1}$ is our small operator in \eqref{small}. Using these medium operators we define $$\calM^\mu(\g)=\calM^\mu:=\langle 1, M_1,\dots,M_r \rangle_{H^{2*}_\G}\subset \calC^\mu$$ the {\em medium algebra}. 

In \cite[Theorem M]{kirillov0} it is  proved that the medium operators are central in $\calC^\mu$. \cite[Theorem 1.1]{kostant-mf}  and the finite dimensional von Neumann double centralizer theorem imply the following
 \begin{theorem}\label{medium} \begin{enumerate}\item For $x\in \s$ the restriction \eqref{restrict} satisfies $\calM^\mu_x=\varrho_\mu(U(\g_x))$. \item $\calM^\mu=Maps(f:\s\to \End(V^\mu) \mid f(x)\in \varrho_\mu(U(\g_x))\subset$ 
$\End(V^\mu))\subset \calC^\mu$. In particular, $\calM^\mu$ is independent of the choice of generating set of $\C[\g]^\G$.  \item The medium algebra $\calM^\mu = Z(\calC^\mu)$ is the center of the Kirillov algebra.
 \end{enumerate}
 \end{theorem}

\subsection{Limits of weight spaces from common eigenspaces of $\M^\mu$}

Denote the maximal torus $\T=\G_{h+e}\subset \G$ corresponding to the centraliser of the regular semisimple element $h+e$. For dominant weights $\mu,\lambda\in \Lambda^+$ we denote by $V^\mu_\lambda\subset V^\mu$ the $\lambda$-weight space of $\T$ in $V^\mu$. Motivated by Kostant's study \cite{kostant} of the zero weight space $V^\mu_0$ Brylinski \cite{brylinski} introduced a filtration \begin{align} \label{filtration} 0<F_0<\dots <F_p<F_{p+1}<\dots<V^\mu_\lambda\end{align} called the {\em Brylinski-Kostant filtration}. It is defined using our regular nilpotent $e\in \g$    as $$F_p:=\{x\in V^\mu_\lambda: e^{p+1} x = 0\}.$$ 
In turn, Brylinski considers the  {\em $e$-limit} of $V^\mu_\lambda$ as \begin{align}\label{limitweight}{\mathrm{lim}}_e V^\mu_\lambda:=\sum e^p\cdot F_p\subset V^\mu.\end{align} The main result of \cite{brylinski} is that $$\small{\sum_p \dim(F_{p+1}/F_{p}) q^p = q^{-(\lambda,\rho)}\sum_k \dim([\lim_e V^\mu_\lambda])^{h=k} q^{\frac{k}{2}}=m^\mu_\lambda(q).}$$ Here $\rho$ is the half-sum of positive roots, $(,)$ is the basic inner product and $[\lim_e V^\mu_\lambda]^{h=k}$ the $k$-eigenspace of $h$ acting on $\lim_e V^\mu_\lambda$. While \begin{align} \label{qanalogue} m^\mu_\lambda(q)=\sum_{w\in W} \epsilon(w)\calP_q(w(\mu+\rho)-\lambda-\rho)\end{align} is Lusztig's \cite{lusztig} $q$-analogue of weight multiplicity. It is defined using the $q$-analogue of Kostant's partition function: $\prod_{\alpha\in \Delta_+} (1-qe^\alpha)^{-1}=\sum_{\pi\in \Lambda} \calP_q(\pi)e^\pi,$ where $\Delta_+\subset \Lambda$ denotes the set of positive roots.

For $z\in \C^\times$, using the $\C^\times$-action \eqref{action}, let $$h_z:=e+zh=z\cdot (e+h)\in \g$$ a regular semisimple element. Define also the  $\C^\times$-action on the Grassmannian $\Gr(k,V^\mu)$ of $k$-planes in $V^\mu$ by $z\cdot U:=\rho_\mu(H^z)(U)\in \Gr(k,V^\mu)$ for $U\in \Gr(k,V^\mu)$. Then we have the following

\begin{theorem}\label{limit}\begin{enumerate} Let $\lambda\leq \mu\in \Lambda^+$, that is $\lambda$ a dominant weight in $V^\mu$, then we have \item \label{firstlimit} for $z\in \C^\times$ the subspace $z\cdot V^\mu_\lambda\subset V^\mu$ is a weight space for the maximal torus $\G_{h_z}$ and thus a common eigenspace for $\M^\mu_{h_z}=\varrho_{\mu}(U(\g_{h_z}))$,
\item \label{secondlimit} $\lim_e V^\mu_\lambda=\lim_{z\to 0} z\cdot V^\mu_\lambda$, i.e.  Brylinski's limit agrees with an actual limit, 
\item $\lim_e V^\mu_\lambda=\lim_{z\to 0} z\cdot V^\mu_\lambda $ is an eigenspace of $\M^\mu_e=\varrho_{\mu}(U(\g_{e}))$ thus  $\lim_e V^\mu_\lambda\subset (V^\mu)^{\G_e}$ \cite[Proposition 2.6]{brylinski},
\item $\lim_e V^\mu_{\mu_{min}}=(V^\mu)^{\G_e}$ for $\mu_{min}$  the minuscule dominant weight in $V^\mu$   (\cite[Corollary 2.7]{brylinski} for $\mu_{min}=0$).
\end{enumerate}
\end{theorem}

\section{Definition and basic properties of big algebras} 

Replacing the symmetric algebra $S^*(\g)$ with the universal enveloping algebra $U(\g)$, Kirillov  in \cite{kirillov0}  also introduced $$Q^\mu(\g)=Q^\mu:=\left(U(\g)\otimes \End(V^\mu)\right)^\G$$ the {\em quantum Kirillov algebra}, which is an algebra over the center $Z(\g)=U(\g)^\G$ of the enveloping algebra. The universal enveloping algebra $U(\g)$ has a canonical filtration $F_0 U(\g)\subset \dots \subset F_k U(\g)\subset F_{k+1}U(\g)\subset \dots$ such that the associated graded algebra $gr(U(\g))\cong S^*(\g)$.  The Rees construction for the filtered algebra $R= U(\g)$ then yields the graded $\C[\hbar]$-algebra \begin{align} \label{rees} R_\hbar:=\oplus_{i=0}^\infty \hbar^{i} F_i R.\end{align} The so-obtained algebra $U_\hbar(\g)$ interpolates between $U_1(\g)\cong U(\g)$ and $U_0(\g)\cong gr(U(\g))\cong S^*(\g).$ We will also consider the {\em $\hbar$-quantum Kirillov algebra} $$\calQ^\mu_\hbar(\g)=\calQ^\mu_\hbar:=\left(U_\hbar(\g)\otimes \End(V^\mu)\right)^\G,$$ which is naturally a $Z_\hbar(\g):=U_\hbar(\g)^\G$-algebra. It interpolates between  the quantum and classical Kirillov algebras: $\calQ^\mu_1\cong\calQ^\mu(\g)$ over $Z_1(\g)=Z(\g)$  and $\calQ^\mu_0\cong\calC^\mu(\g)$ over $S^*(\g)^\G\cong Z_0(\g)$.

Recall from \cite{feigin-frenkel-reshetikhin, rybnikov} and specifically from \cite[\S 8.2]{yakimova} the two-point Gaudin algebra $\calG\subset\calQ(\g):=(U(\g)\otimes U(\g))^\G.$ This is defined as a quotient of the Feigin-Frenkel center \cite{feigin-frenkel}, and thus it is a commutative subalgebra of the universal quantum Kirillov algebra $\calQ(\g)$. We will also take the Rees construction \eqref{rees} with respect to the filtration on $\calQ$ and $\calG$ coming from the filtration on the first copy of $U(\g)$ and denote them $\calG_\hbar\subset \calQ_\hbar$. These are graded $\C[\hbar]$-algebras, with  $\calG_\hbar$ commutative. For $\mu\in \Lambda^+(\G)$ the image $\calG_\hbar^\mu:=\pi^\mu(\calG_\hbar)\subset \calQ^\mu_\hbar$ under the projection $\pi^\mu:\calQ_\hbar\to \calQ^\mu_\hbar$ induced from the projection $U(\g)\to \End(V^\mu)$ is called the {\em $\hbar$-quantum big algebra}, which interpolates between $\calG^\mu:=\calG^\mu_1\subset \calQ^\mu$ the {\em quantum big algebra} and $\calB^\mu:=\calG^\mu_0\subset \calC^\mu$ the {\em (classical) big algebra}. 

The universal big algebras $(\calG_0)_x$ for $x\in \s$ were denoted by $\calA_x\subset U(\g)$ in \cite{feigin-frenkel-rybnikov} and its action on a representation $V^\mu$ was also studied in {\em loc. cit.}. Our  finite dimensional matrix algebras $\calB^\mu_x$  from \eqref{restrict} are just the images of $\calA_x$ in $\End(V^\mu)$. Using their results we can deduce the following

\begin{theorem} \label{big} Let $\mu\in \Lambda^+(\G)$ be a dominant character. Then \begin{enumerate} \item the $\hbar$-quantum big algebra $\calG_\hbar^\mu\subset \calQ^\mu_\hbar$ is a maximal commutative subalgebra, finite-free over $Z_\hbar(\g)$, consequently it contains the {\em $\hbar$-quantum medium algebra} $\M_\hbar^\mu:=Z(\calQ^\mu_\hbar)\subset \calG^\mu_\hbar$, 
\item the  big algebra $\calB^\mu=\calG^\mu_0\subset \calC^\mu$ is a maximal commutative subalgebra, finite-free over $S^*(\g)^\G$, consequently, the medium algebra $\M^\mu\cong \M^\mu_0\cong Z(\calC^\mu)\subset \calB^\mu$, \item the Hilbert series of $\calB^\mu$ satisfies $$\sum_{i=0}^\infty \dim((\calB^\mu)^i)q^i=\frac{\prod_{\alpha\in \Delta^+}\frac{(1-q^{(\rho+\mu,\alpha)})}{(1-q^{(\rho,\alpha)})}}{\prod_{j=1}^r (1-q^{d_j})},$$
\item \label{thirdbig} for all $x\in \s$ the algebra $\calB_x^\mu\subset  \End(V^\mu)$ acts both with $1$-dimensional common eigenspaces and cyclically. 
\end{enumerate}
\end{theorem}
It was already observed in \cite{feigin-frenkel-rybnikov} that Theorem~\ref{big}.\ref{thirdbig}  implies that the cyclic action of $\calB^\mu_e$ on $V^\mu$ endows $V^\mu$ with a graded ring structure. The whole big algebra $\calB^\mu$ however contains much more information. For example it follows from Theorem~\ref{limit}.\ref{firstlimit} that $\calB^\mu_{h_z}$ leaves $z\cdot V^\mu_\lambda$, the common eigenspaces of $\M^\mu_{h_z}=\varrho(U(\g_{h_z}))\subset \End(V^\mu)^{\G_{h_z}}$, invariant. Thus by Theorem~\ref{limit}.\ref{secondlimit}  $\calB^\mu_e$ leaves $\lim_e V^\mu_\lambda$ invariant and so we can define the {\em multiplicity algebra} \begin{align}\label{multiplictyalgebra}Q^\mu_\lambda:=\calB^\mu_e\mid_{\lim_e V^\mu_\lambda}\subset \End(\lim_e V^\mu_\lambda).\end{align}
Then Theorem~\ref{big}.\ref{thirdbig} and Theorem~\ref{limit} imply the following

\begin{corollary} \label{multiplicity} Let $\lambda\leq \mu\in \Lambda^+(\G)$ be dominant characters. The big algebra $\calB^\mu_e$ at $e\in \s$ induces \eqref{multiplictyalgebra} a graded algebra structure $Q^\mu_\lambda$ on $(\lim_e V^\mu_\lambda)^*$ such that \begin{enumerate} \item $\sum\dim(Q^\mu_\lambda)^iq^{(\mu-\lambda,\rho)-i}=m^\mu_\lambda(q)$ Lusztig's $q$-analogue of multiplicity \eqref{qanalogue}, \item there are natural quotient maps $\calB^\mu_e\twoheadrightarrow Q^\mu_{\mu_{min}}\twoheadrightarrow Q^\mu_\lambda$,  \item \label{thirdmultiplicity}  $Q^\mu_{\mu_{min}}\cong \calB_e^\mu/((\calM^\mu_{e})_+)=\calB_e^\mu/((M_1)_e,\dots,({M_r})_e)$.
\end{enumerate}
\end{corollary}

\subsection{Computing big algebras}

Fix a basis $\{X_i\}$ for $\g$ and a dual basis $\{X^i\}\subset \g$ with respect to the Killing form of $\g$. For $A\in \calC^\mu$, following Kirillov \cite{kirillov0}, Wei \cite{wei} introduced the following $D$-operator: 
$$D(A):=\frac{1}{2}\sum_i \rho_\mu(X^i){\frac{\partial(A)} {\partial X_i}}.$$ It is shown in \cite{wei} that $D(A)\in \calC^\mu$ and that $D(A)$ is independent of the choice of the basis $\{X_i\}\subset \g$.  This $D$-operator allows us to construct  new operators from known ones. For example for  $p\in \C[\g]^\G$  we have   $D(p)=M_{dp/2}$ is the medium operator of \eqref{mediumdef}.   It is not true that for any $p\in \C[\g]^\G$ iterated derivatives $D^k(p)$ are still in the big algebra $\calB^\mu$. However starting with a good generating set of $\C[\g]^\G$ we can explicitly generate the big algebra.
    Here is such an example in type $A$.

    \begin{theorem}\label{explicit} For $A\in \sl_n$ let $c_i(A)=(-1)^i(\det(\Lambda^i(A))$ be the $i$th coefficient of the characteristic polynomial of $A$. Then $\C[\sl_n]^{\SL_n}\cong \C[c_2,\dots,c_n]$ and the {\em big operators} $$B_{i,k-i}=D^i(c_k)\in \calC^\mu$$ generate the big algebra $$\calB^\mu=\C[\sl_n]^{\SL_n}\langle B_{i,k-i}\rangle_{0<i<k\leq n}\subset \calC^\mu.$$
    \end{theorem}

    Similar generating sets are known in types $B,C,D,G$ and conjectured to exist in all types \cite{yakimova}.

\section{Geometric aspects}

Let $\G$ be a connected semisimple complex Lie group, $\G^\vee$ its Langlands dual group.  Their Lie algebras are $\g$ and $\g^\vee$ and $\t\subset \g$ and $\t^\vee\subset \g^\vee$ are Cartan subalgebras with $\t^*\cong \t^\vee$ naturally.  Identify $\g\cong\g^*$ and $\t\cong \t^*$ by the Killing form. Then  the Duflo isomorphism \cite[Lemme V.1]{duflo} is\begin{align} \label{duflo}\delta:=\chi^{-1}\circ \psi: Z(\g)\to S^*(\t)^\W\cong   S^*(\g)^\G \end{align} where $\chi: S^*(\g)^\G\cong \C[\g]^\G\to \C[\t]^\W\cong S^*(\t)^\W$ is the Chevalley isomorphism and $\psi:Z(\g)\to S^*(\t)^\W$ is the Harish-Chandra isomorphism. On the Rees constructions \eqref{rees} this induces
\begin{align*}\delta_\hbar:Z_\hbar(\g)\cong \C[\hbar][\g]^\G\cong \C[\g^\vee\times \C]^{\G^\vee \times \C^\times}.\end{align*}

The following Theorem~\ref{geometric} shows that  our algebras have natural meanings related to  equivariant (intersection) cohomology of affine Schubert varieties. All our cohomologies and intersection cohomologies will be with $\C$-coefficients and $\G$-equivariant (intersection) cohomology will be over $H^{2*}(B\G)\cong \C[\g]^\G=H^{2*}_\G$. From results in \cite{bezrukavnikov-finkelberg} we can deduce the following

\begin{theorem}\label{geometric} Let $\G$ be a connected semisimple group and $\g$ its Lie algebra, with Langlands dual $\G^\vee$ and corresponding affine Grassmannian $\Gr:=\Gr_{\G^\vee}=\G^\vee(\C((z)))/\G^\vee(\C[[z]])$. Let $\mu\in \Lambda^+(\G)$ be a dominant character and let $\Gr^\mu:=\overline{\G^\vee(\C[[z]])z^\mu}$ be the corresponding {\em affine Schubert variety}, with action of  $\G^\vee\subset \G^\vee(\C[[z]])$ from the left and $\C^\times$ through loop rotation on $z$.  For $\lambda\leq \mu \in \Lambda^+(\G)$ we let  $\mathcal{W}^\mu_\lambda:=\G^\vee_1(\C[[z^{-1}]])z^\lambda\cap \Gr^\mu$ be the {\em affine Grassmannian slice}, where $\G^\vee_1(\C[[z^{-1}]])$ is the kernel of the evaluation map $\G^\vee(\C[[z^{-1}]])\to \G^\vee$ at $z^{-1}=0$. Then \begin{enumerate} \item $H^{2*}_{\G^\vee \times \C^\times}(\Gr^\mu)\cong \calM_\hbar^\mu$ as $H^{2*}_{\G^\vee \times\C^\times}\cong \C[\hbar][\g]^\G$-algebras, \item \label{nakajima}
 $\End_{H^{2*}_{\G^\vee \times \C^\times}(\Gr^\mu)}(I\!H^{2*}_{\G^\vee \times \C^\times}(\Gr^\mu))\cong \calQ_\hbar^\mu$ as $\H^{2*}_{\G^\vee \times \C^\times} \cong Z_\hbar(\g)$-algebras, \item  $I\!H^{2*}_{\G^\vee \times \C^\times}(\Gr^\mu)\cong \calG_\hbar^\mu$ as  $H_{\G^\vee \times \C^\times}^*(\Gr^\mu)\cong\calM_\hbar^\mu$-modules. In particular, $\calG^\mu_\hbar$ endows $I\!H^{2*}_{\G^\vee \times \C^\times}(\Gr^\mu)$ with a graded ring structure compatible with the action of $H_{\G^\vee\times \C^\times}^{2*}(\Gr^\mu)\cong\calM_\hbar^\mu$,
 \item $I\!H^{2*}(\mathcal{W}^\mu_\lambda)\cong Q^\mu_\lambda$ as graded vector spaces, thus $Q^\mu_\lambda$ endows $I\!H^{2*}(\mathcal{W}^\mu_\lambda)$
 with a graded ring structure. 
 \end{enumerate}
\end{theorem}

\section{Examples--problems}
\subsection{Minuscule and weight multiplicity free Kirillov algebras} When $V^\mu$ is weight multiplicity free, for example when $\mu$ is minuscule,  the Kirillov algebras are already commutative \cite[Theorem 4.1]{rozhkovskaya}, thus $\calM_\hbar^\mu\cong\calG_\hbar^\mu\cong \calC_\hbar^\mu$. First we discuss the classical case of $\calB^\mu=\calG^\mu_0$. 

For any $\mu\in \Lambda^+$  we have the unique closed $\G$-orbit $\G v_\mu\cong \G/\rP_\mu\subset \mathbb{P}(V^\mu)$, a partial flag variety. We can form the {\em big zero scheme} $\calZ^\mu:=\cap_{B\in \calB^\mu} \calZ(Y_B)\subset \s \times \mathbb{P}(V^\mu)$ 
as the common zeroes of the vector fields $Y_B\in \mathfrak{X}(\s \times \mathbb{P}(V^\mu))$ induced by  the big operators $B\in\calB^\mu$, parametrizing their common eigenvectors. By construction $\C[\calZ^\mu]\cong \calB^\mu$. On the other hand we can see  that $\calZ(Y_{M_1})\cap \G v_\mu\subset \calZ^\mu\subset \s\times \mathbb{P}(V^\mu)$, because  for a generic $x\in \s$ the scheme $\calZ ((Y_{M_1})_x)\cap \G v_\mu$ contains only isolated points of $\calZ ((Y_{M_1})_x)$.  From  \cite[Theorem 1.3]{hausel-rychlewicz} we have that $\C[\calZ(Y_{M_1})\cap \G v_\mu]\cong H^{2*}_\G(\G v_\mu)$ and thus we always have a surjective map \begin{align}\label{ghc} \calB^\mu\twoheadrightarrow H^{2*}_\G(\G/\rP_\mu).\end{align} The ring homorphism \eqref{ghc} can be thought of  an upgrade of a similar linear map $\hat{f}$ in \cite[Theorem 1]{friedman-morgan}, which was proved (essentially) in \cite{ginzburg} to be a surjection.   When $\mu$ is minuscule, the Hilbert series of the two graded rings of \eqref{ghc} agree and we get that 
$ \calB^\mu\cong H^{2*}_\G(\G/\rP_\mu).$  
This result was deduced by algebraic means in \cite[\S 6]{panyushev}.

When we use the $\hbar=0$ specialization of\ Theorem~\ref{geometric}.1 we get that \begin{align} \label{cominuscule} \calB^\mu\cong \calM^\mu \cong H_{\G^\vee}^{2*}(\Gr^\mu)\cong  H_{\G^\vee}^{2*}(\G^\vee/\rP^\vee_\mu)\end{align} the equivariant cohomology of the cominuscule flag variety. The two descriptions above then agree because $\H^{2*}_\G(\G/\rP_\mu)\cong \C[\t]^{\W_\mu}\cong \C[\t^*]^{\W_\mu}\cong \C[\t^\vee]^{\W_\mu}\cong  \H^{2*}_{\G^\vee}(\G^\vee/P^\vee_\mu),$ where $\W_\mu:=Stab(\mu)\subset \W$ in the Weyl group  of $\G$.

Similarly, for $V^\mu$ weight multiplicity free \cite[Conjecture 6]{panyushev} suggests $\G$-invariant subvarieties $X_\mu\subset \mathbb{P}(V^\mu)$ such that $\calB^\mu\cong H^{2*}_\G(X_\mu)$. For example for the weight multiplicity free $\mu=k\omega_1\in \Lambda^+(\SL_n)$ we have
 $X_\mu\cong S^k(\P^{n-1})$,  the $k$th symmetric product with the diagonal action of $\SL_n$. With a similar technique as above and straightforwardly extending
 \cite[Theorem 1.3]{hausel-rychlewicz} to the orbifold $S^k(\P^{n-1})$ we can prove Panyushev's conjecture:
 \begin{align}\label{symmetric}\calB^{k\omega_1}(\sl_n)\cong H^{2*}_{\SL_n}(S^k(\P^{n-1}))\cong S^k_{H^{2*}_{\SL_n}}(H^{2*}_{\SL_n}(\P^{n-1})).\end{align}
Note that $\calB^{k\omega_1}(\mathfrak{s}\mathfrak{l}_n)\cong H^{2*}_{\PGL_n}(\Gr^{k\omega_1})$ from Theorem~\ref{geometric}.\ref{nakajima}. The varieties $\Gr^{k\omega_1}$ are different from $S^k(\P^{n-1})$
for example   $S^k(\P^{1})\cong \P^k$ is smooth while $\Gr^{k\omega_1}(\PGL_2)$ is  singular for $k>1$. Still they have isomorphic equivariant cohomology rings: \begin{align}\label{sl2} H^{2*}_{\SL_2}(\P^k)\cong \calB^{k\omega_1}(\mathfrak{s}\mathfrak{l}_2)\cong H^{2*}_{\PGL_2}(\Gr^{k\omega_1}).\end{align}


For quantum Kirillov algebras Theorem~\ref{geometric}.\ref{nakajima} is useful when $
\mu$  is minuscule. In that case the loop rotation action on $\Gr^\mu$ is trivial, which implies the surprising
\begin{corollary}\label{minuscule} When $\mu\in \Lambda^+(\G)$ is minuscule $\calC^\mu(\g)\cong\calQ^\mu(\g)$ as $Z(\g)\stackrel{\delta}{\cong} \C[\g]^\G$-algebras, where $\delta$ of \eqref{duflo} is the Duflo isomorphism. \end{corollary} 

The isomorphism can be constructed as the combination of the generalised Harish-Chandra isomorphisms in \cite[\S 9]{higson}, making it the sought-after generalised Duflo isomorphism in this minuscule case.

Applied to the standard representation $\calQ^{\omega_1}(\sl_n)\cong\calC^{\omega_1}(\sl_n)$ Corollary~\ref{minuscule} implies that the Capelli identity matches the classical Cayley-Hamilton identity under the Duflo isomorphism, which is \cite[Theorem 7.1.1]{molev}.  In types $C$ and $D$  the case of $N=2n$ in \cite[ Theorem 7.1.6]{molev} gives  $\calQ^{\omega_1}(\g)
\cong \calC^{\omega_1}(\g)$. Note that in type $B$, the standard representation is not minuscule. Indeed the  case of $N=2n+1$ in \cite[Theorem 7.1.6]{molev} shows that the quantum Capelli identity does not map to the classical Cayley-Hamilton equation, thus $\calQ^{\omega_1}(\so_{2n+1})\ncong \calC^{\omega_1}(\so_{2n+1})$, which is compatible with the non-triviality of the loop rotation on $\Gr^{\omega_1}(\SO_{2n+1})$.

\subsection{Visualisation of explicit examples}
\label{explicitsecion}

 As the big algebras $\calB^\mu$  are commutative  and finite-free over the polynomial ring $H^{2*}_\G$, they correspond to affine schemes $\Spec(\calB^\mu)$ finite flat over the affine space $\Spec(H^{2*}_\G)$. With the exception of some small rank examples the embedding dimension of $\Spec(\calB^\mu)$ (the minimal number of generators of $\calB^\mu$) is larger than three, thus we cannot directly depict them. For visualisation purposes the {\em principal subalgebras} obtained by base changing to a principal $\SL_2\to \G$ subgroup: $\calB^\mu_{\SL_2}:=\calB^\mu\otimes_{H^{2*}_\G} \H^{2*}_{\SL_2}$ and $\calM^\mu_{\SL_2}:=\calM^\mu\otimes_{H^{2*}_{\G}} \H^{2*}_{\SL_2}$ are better behaved. Their spectra 
$\Spec(\calB^\mu_{\SL_2})$ and $\Spec(\calM^\mu_{\SL_2})$, which we call the {\em big and medium skeletons}, are curves over the line $ \Spec(H^{2*}_{\SL_2})$. We call $\spec(\calB^\mu_{h})$ and $\Spec(\calM^\mu_{h})$,  the fibers over the principal semisimple element $h\in \sl_2/\! /\SL_2\cong \Spec(H^{2*}_{\SL_2})$, the {\em big and medium principal spectra}. Because of Theorem~\ref{medium} one can identify \begin{align}\label{weights}
\Spec(\calM^\mu_{h})\cong \Spec(V^\mu)\subset \t^*
\end{align}
where $\Spec(V^\mu)$ is the reduced scheme of the set of weights in $V^\mu$, which appeared in a closely related context in \cite[Theorem 1.3.2]{ginzburg}.

\subsubsection{Big algebras for $\SL_2$}
By \eqref{sl2} we have $ \calB^{n\omega_1}(\sl_2)\cong H^{2*}_{\SL_2}(\P^{n})$, which have been computed in \cite[\S 4.4]{hausel-rychlewicz}, yielding
\begin{align*}
 \calB^{n\omega_1}(\sl_2)\cong
 \begin{cases}
  \C[c_2,M_1]/\big((M_1^2+n^2c_2)(M_1^2+(n-2)^2c_2)\dots(M_1^2+4c_2)M_1\big) &\text{ for $n$ even;}\\
  \C[c_2,M_1]/\big((M_1^2+n^2c_2)(M_1^2+(n-2)^2c_2)\dots(M_1^2+9c_2)(M_1^2+c_2)\big) &\text{ for $n$ odd.}
 \end{cases}
 \end{align*}

\begin{figure}[ht!]
\begin{center}
 {
  \includegraphics[width=7cm]{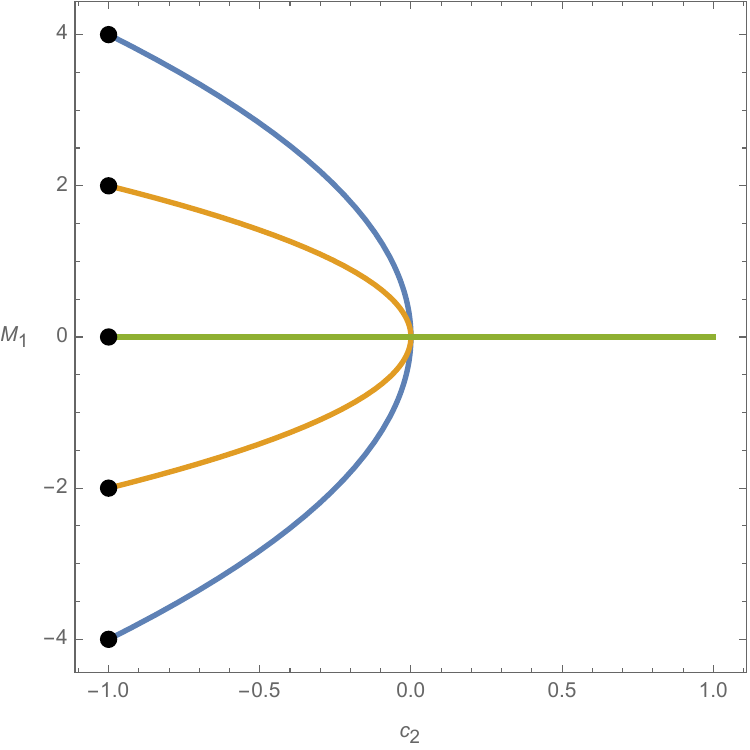}}
   \hskip1.3cm
{
  \includegraphics[width=7cm]{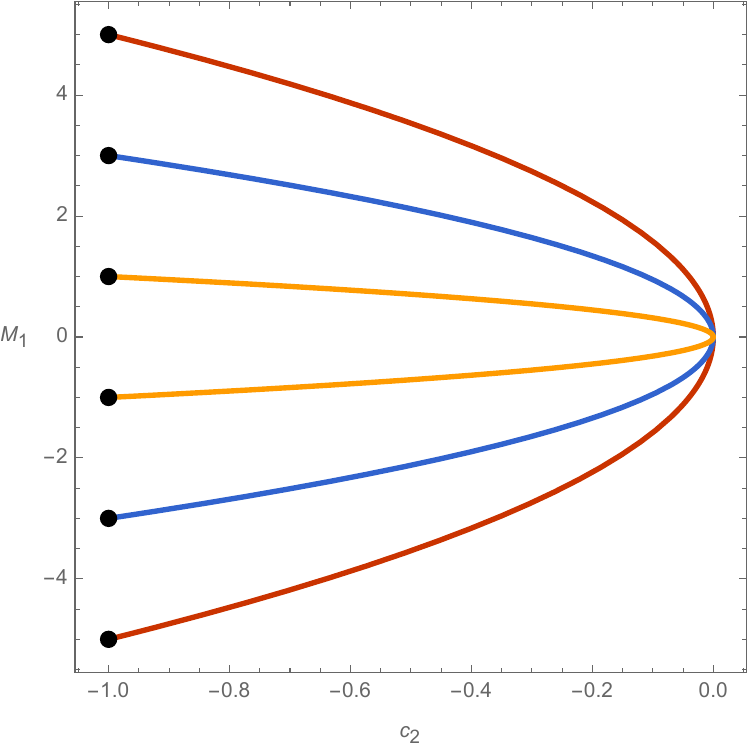}}
\end{center}
\caption{$\Spec \calB^{4\omega_1}(\sl_2)\cong \Spec H_{\SL_2(\C)}^*(\PP^4)$ and $\Spec \calB^{5\omega_1}(\sl_2)\cong\Spec H_{\SL_2(\C)}^*(\PP^5)$.}
\label{sl2pn}
\end{figure}

In Figure~\ref{sl2pn} the real points of the spectrum of the big algebras for two $\SL_2$ examples are shown, with the black dots depicting the principal spectrum, which by \eqref{weights} can be identified with the weights of the representation.

\subsubsection{Big algebra for standard representation of $\SL_3$}

\begin{figure}[ht!] 
\begin{center}
  \includegraphics[width=6cm]{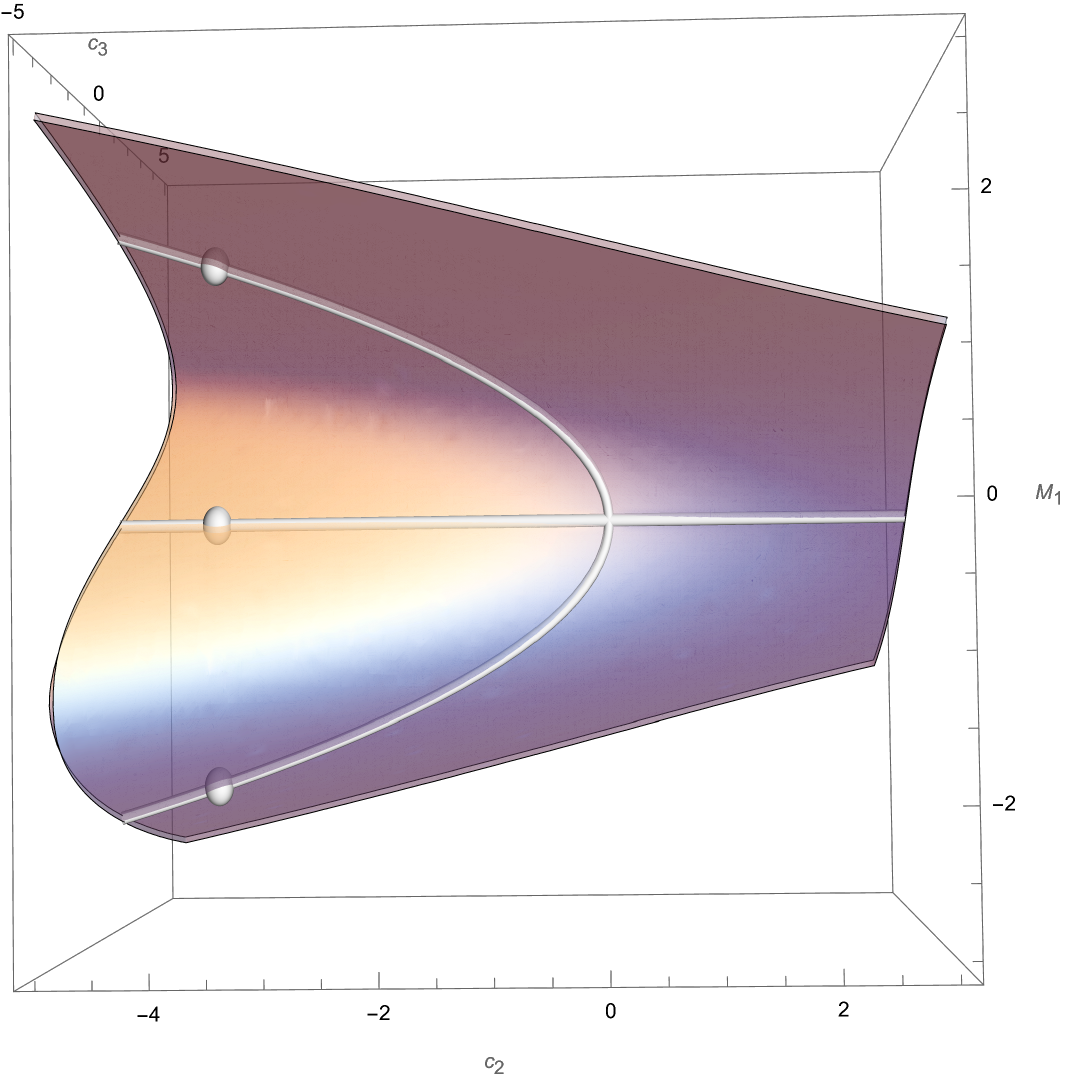}
\end{center}
\caption{$\Spec \calB^{\omega_1}(\sl_3)$, its skeleton $\Spec \calB_{\SL_2}^{\omega_1}(\sl_3)$ and principal spectrum $ \Spec \calB^{\omega_1}_{h}(\sl_3)$.}
\label{tripletskeleton}
\end{figure}

Using the Cayley-Hamilton identity one can explicitly  compute the big algebra for the standard representation of $\SL_3$ in terms of the small operator $M_1$ of \eqref{small} as $$\calB^{\omega_1}(\sl_3)\cong \C[c_2,c_3,M_1]/(M_1^3+c_2 M_1+c_3).$$
Figure~\ref{tripletskeleton} shows the real points of the spectrum of $\calB^{\omega_1}(\sl_3)$ together with its skeleton and principal spectrum. 

\subsubsection{Big algebra for $\rho_{3\omega_1}$ of $\SL_3$ -- the decuplet}

\begin{figure}[ht!]
\begin{center}
{
  \includegraphics[width=5.3cm]{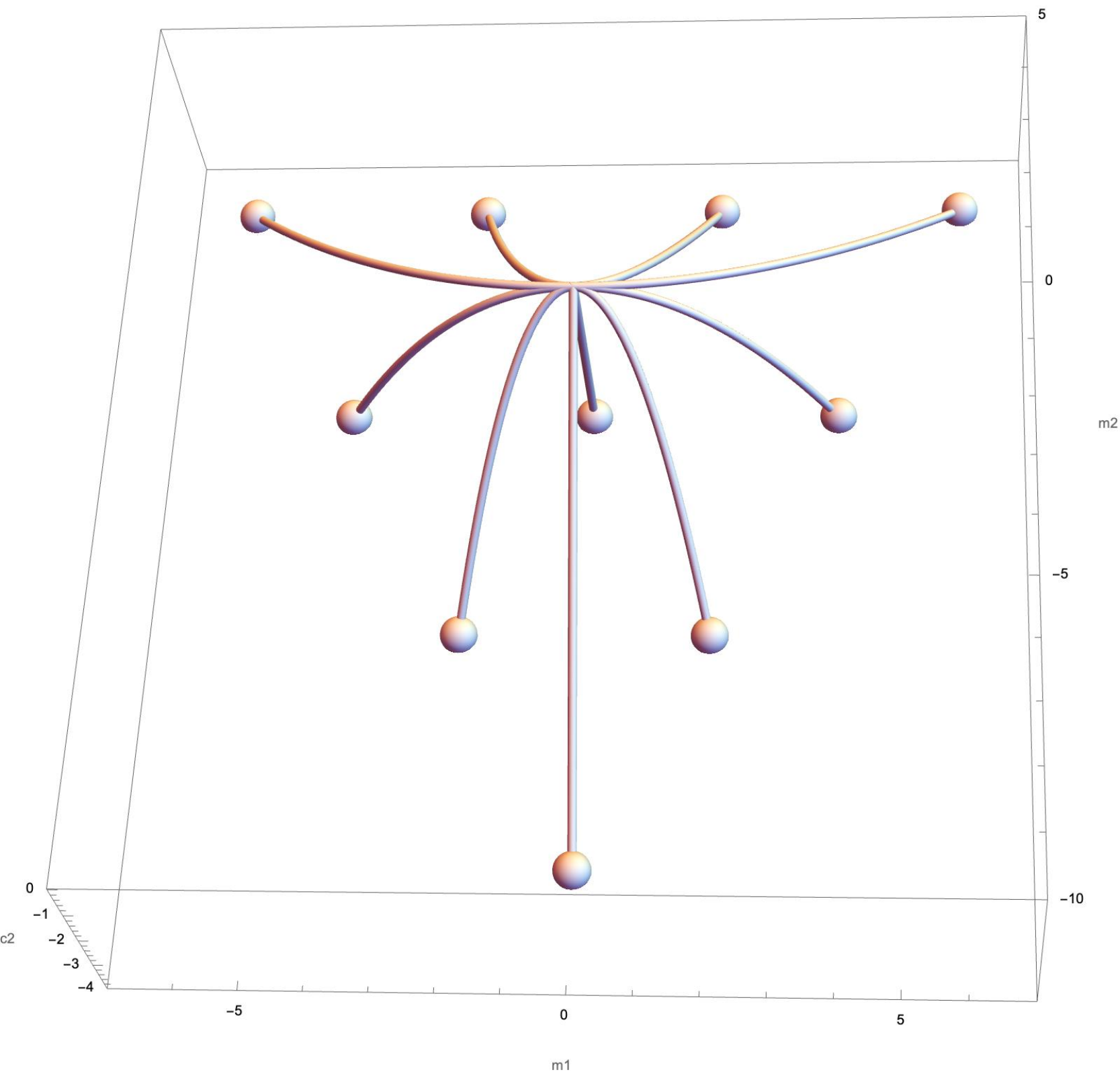}}
  \hfill
{
\includegraphics[width=5.5cm]{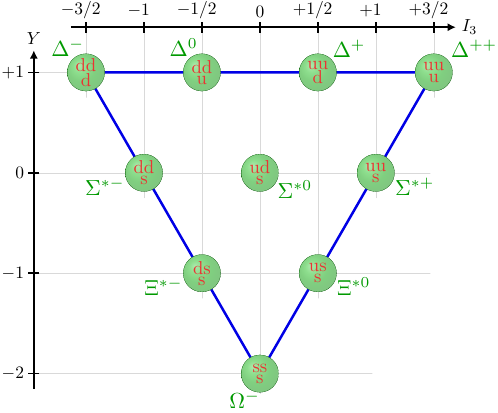}}
   \hfill
{
  \includegraphics[width=5.3cm]{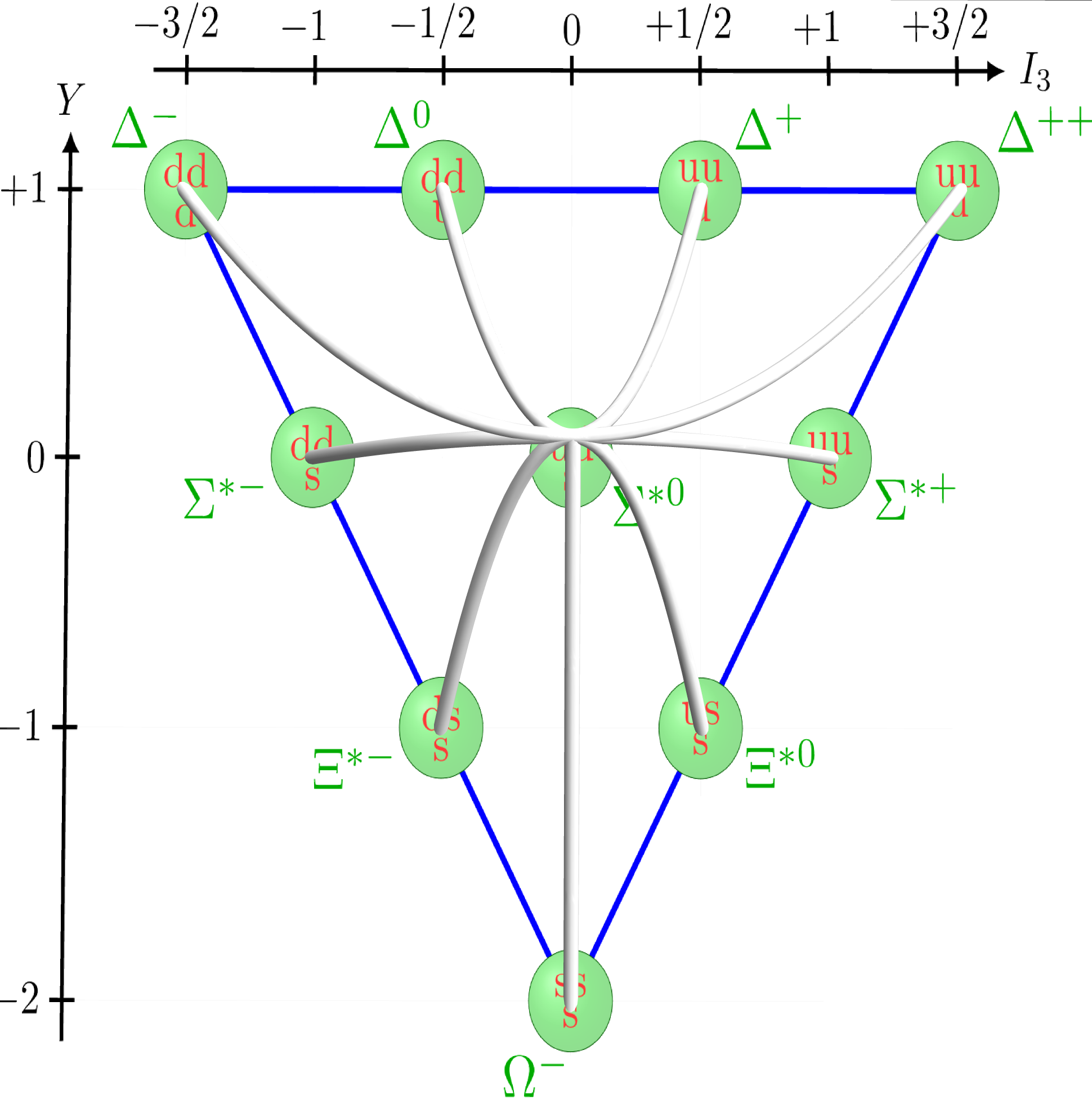}}
\end{center}
\caption{$\Spec(\calB_{\SL_2}^{3\omega_1}(\sl_3))$ over $\Spec(\calB_{h}^{3\omega_1}(\sl_3))$,  baryon decuplet and skeleton over  decuplet}
\label{decuplet}
\end{figure}

Using either \eqref{symmetric} or Theorem~\ref{explicit} we can compute the big algebra $\calB^{3\omega_1}(\sl_3)\cong H^{2*}_{\SL_3}(S^3(\P^{2}))$ explicitly in terms of the medium operators $M_1=D(c_2)$ and $M_2=D(c_3)$:
\begin{align} \label{decupletideal} & \calB^{3\omega_1}(\sl_3)\cong \C[c_2,c_3,M_1,M_2]/\nonumber\\ &
\left({ \begin{array}{l}  M_1^4 - 6M_1^2M_2 + 4M_1^2c_2 - 18M_1c_3 + 3M_2^2 - 6M_2 c_2, \\ M_1^3M_2 + M_1^3c_2 + 3M_1^2c_3  - 3M_1M_2^2+\\ + M_1M_2c_2 + 4M_1c_2^2 - 9M_2c_3 
\end{array}}\right)
\end{align}
 From this we  obtain $\calB_{\SL_2}^{3\omega_1}$ by setting $c_3=0$ and $\calB_{h}^{3\omega_1}$ by further setting $c_2=-4$. The first picture of Figure~\ref{decuplet}  shows the resulting picture of the real points of the skeleton and the principal spectrum.

The principal spectrum can be identified with the set of weights in $V^{3\omega_1}$ by \eqref{weights}, which in turn corresponds to the particles appearing in the baryon decuplet of Gell-Mann \cite[pp. 87, Fig. 1 pp.88]{gell-mann-neeman}, see the second picture in Figure~\ref{decuplet}. There are two quantum numbers, the isospin $I_3$ and hypercharge $Y$ which distinguish  the particles in the multiplet. They correspond to our operators as $(M_1)_h=4I_3$ and $(M_2)_h=4Y$. Thus our two relations  in our big algebra \eqref{decupletideal} give the following  generating set of polynomial relationships between these two quantum numbers in the baryon decuplet:
\begin{align}\label{decupletid}\begin{array}{c}  I_3(Y - 1)(4I_3^2 - 3Y - 4) =0 \\ 16I_3^4 - 24I_3^2Y - 16I_3^2 + 3Y^2 + 6Y =0 \end{array}\end{align}
The third picture in Figure~\ref{decuplet} shows that we can obtain the skeleton $\spec(\calB^{3\omega_1}_{\SL_2})$ by connecting the particles in the decuplet by parabolas when they correspond to each other under the up-down quark symmetry. The two particles fixed by this symmetry, the $\Sigma^{*0}$ and $\Omega^-$, are supporting lines  in the skeleton $\Spec(\calB^{3\omega_1}(\sl_3))$. $\Omega^-$ is the  particle formed by three strange quarks, whose existence was famously predicted by Gell-Mann  based on this  baryon decuplet model \cite[pp. 87]{gell-mann-neeman}. 

\subsubsection{Big algebra of adjoint representation of $\SL_3$ -- the octet}

The smallest dimensional non-weight multiplicity free representation is the adjoint representation $\rho_{\omega_1+\omega_2}$ of $\SL_3$. In this case $\calM^{\omega_1+\omega_2} (\sl_3)\subsetneq \calB^{\omega_1+\omega_2} (\sl_3)$, the medium and big algebras are distinct. Using \cite[Table III]{rozhkovskaya}  or Theorem~\ref{explicit} one can compute the big algebra, and in turn the medium subalgebra,  explicitly, in terms of the medium operators $M_1=D(c_2)$ and $M_2=D(c_3)$ and big operator $N_1=D^2(c_3)$:   
\begin{align}\label{octetbigideal}
& \calB^{\omega_1+\omega_2}(\sl_3)\cong \C[c_2,c_3,M_1,N_1]/\nonumber \\ &
\left(\begin{array}{l} 3M_1^2+N_1^2+12c_2, \\  M_1^3 N_1+c_2M_1N_1-9 c_3M_1
\end{array}\right)\\ \label{octetmedideal} & \calM^{\omega_1+\omega_2}(\sl_3)\cong \C[c_2,c_3,M_1,M_2]/ \nonumber\\ &
\left(\begin{array}{l} M_1^2M_2+c_2M_2+3c_3M_1, \\M_1^4+4c_2M_1^2+3M_2^2, \\
3M_1M_2^2+9c_3M_2-c_2M_1^3-4c_2^2M_1
\end{array}\right)
\end{align}

\begin{figure}[bt!] 
\begin{center}
 {
  \includegraphics[width=5.3cm]{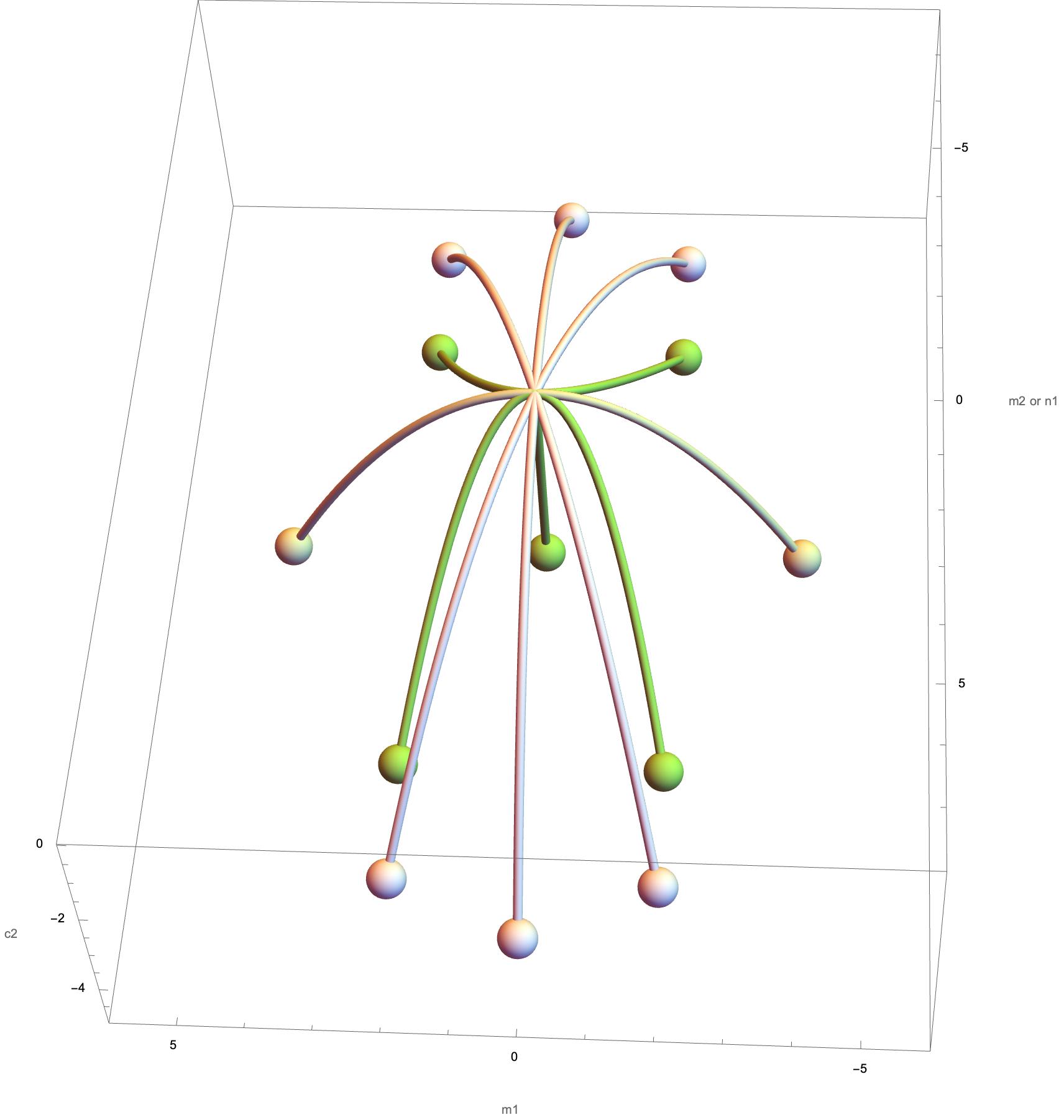}}
  \hfill
{
  \includegraphics[width=5.5cm,trim=0cm 0cm 0 0]{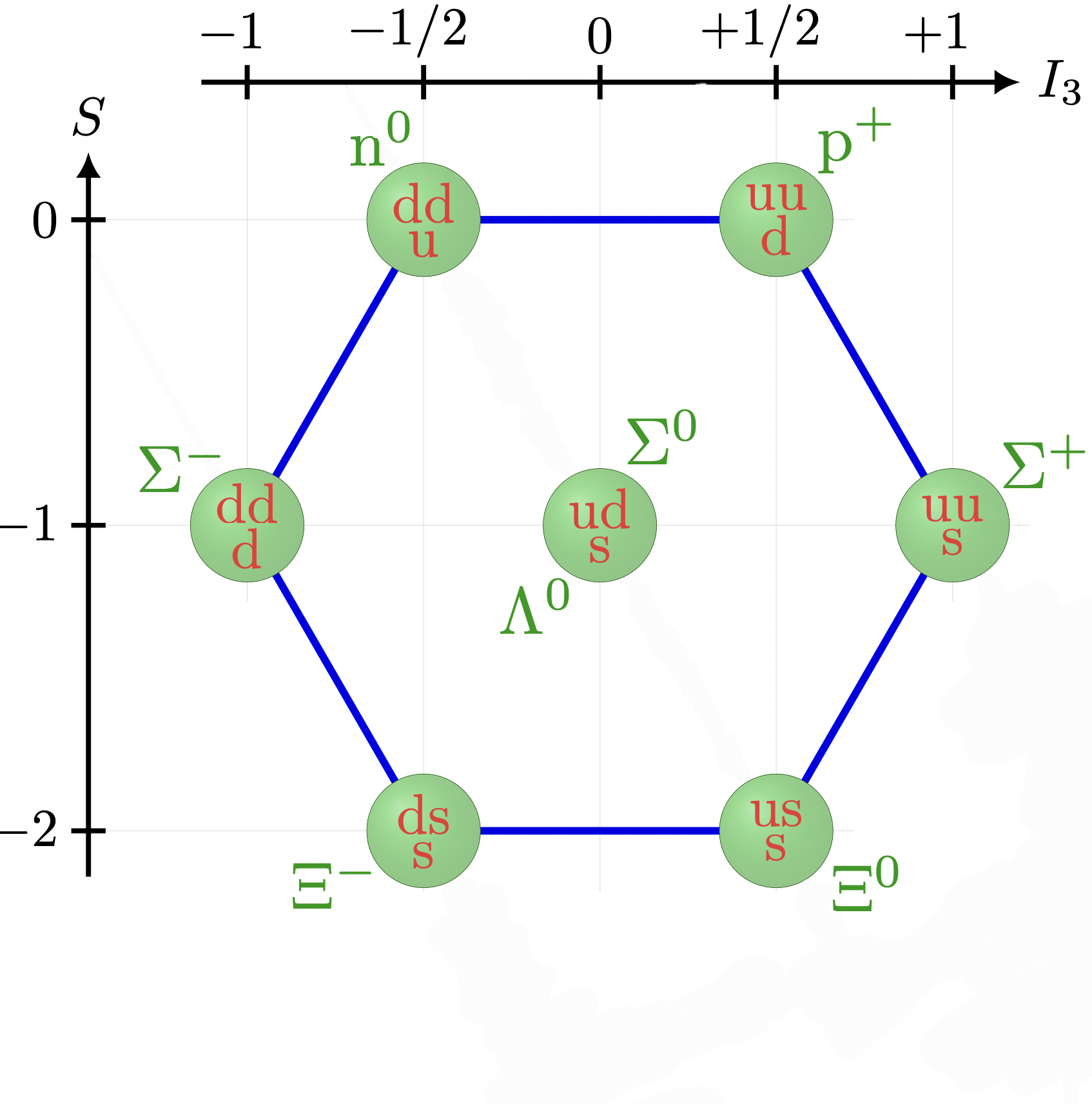}}
    \hfill
{
\includegraphics[width=5.3cm,trim=0cm 0cm 0cm 0]{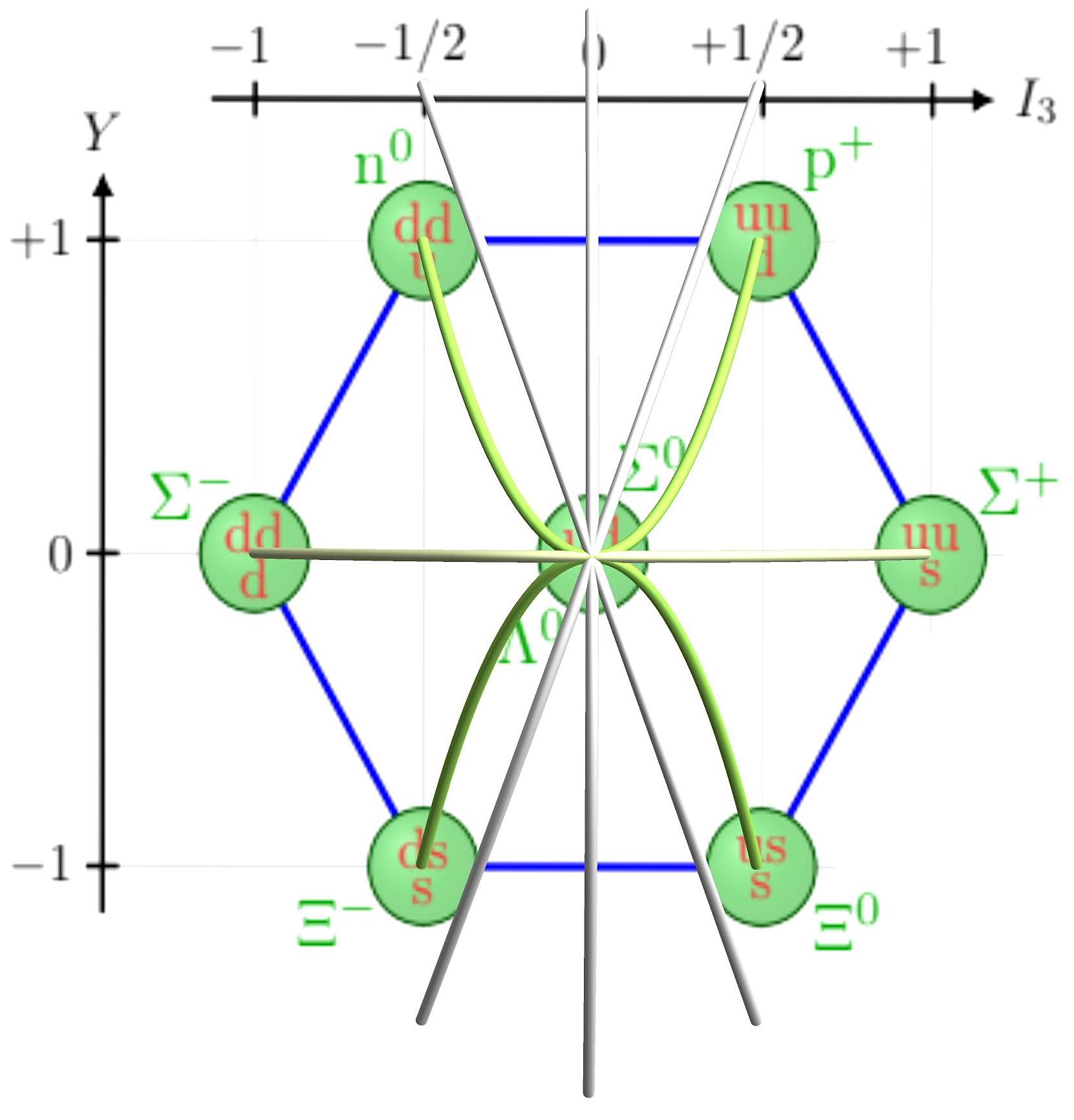}}
\end{center}
\caption{Skeletons $\calB_{\SL_2}^{\omega_1+\omega_2}(\sl_3)$, $\calM_{\SL_2}^{\omega_1+\omega_2}(\sl_3)$ over $\calB_{h}^{\omega_1+\omega_2}(\sl_3)$, $\calM_{h}^{\omega_1+\omega_2}(\sl_3)$,  baryon octet and big and medium skeletons over  octet}
\label{octetfigure}
\end{figure}

Setting $c_3=0$ in these equations gives us the big and medium skeletons, why further specialising $c_2=-4$ gives us the big and medium principal spectra. These are depicted (white for big and green for medium) on the first picture of Figure~\ref{octetfigure}. We used the coordinates $c_2,M_1$ and $N_1$ for the big skeleton but $c_2,M_1$ and $M_2=\frac{1}{3} M_1N_1$ for the medium skeleton. 

Thus our relations in \eqref{octetmedideal} imply the following generating set of polynomial relations between the quantum numbers $I_3$ and $Y$ in the baryon octet (see second picture in Figure~\ref{octetfigure}):

\begin{align}\label{octetid}\begin{array}{c} 
    Y(2I_3 - 1)(2I_3 + 1) =0\\4I_3^3 + 3I_3Y^2 - 4I_3=0\\ 16I_3^4 - 16I_3^2 + 3Y^2=0 \end{array}
    \end{align}

    We can also compute the multiplicity algebra of the $0$ weight from \eqref{octetbigideal} and Corollary~\ref{multiplicity} to get $$Q_0^{\omega_1+\omega_2}(\sl_3)\cong \calB^{\omega_1+\omega_2}/((\calM^{\omega_1+\omega_2})_+)\cong \C[N_1]/(N_1^2).$$

On the third picture of Figure~\ref{octetfigure} we can see that the medium skeleton can be built on the baryon octet by connecting the particles corresponding by up-down quark symmetry -- such as the neutron $n^0$ and proton $p^+$ -- with parabolas. The big skeleton is more complicated. It consists of four parabolas (one  shared with the medium skeleton) and has two points in its principal spectrum over the origin in the baryon octet corresponding to the multiplicity two $0$ weight space containing the two particles $\Sigma^0$ and $\Lambda^0$. 

\begin{remark} Using  \cite{tai}, where the Kirillov algebra is computed for the adjoint representation of any simple complex Lie group, one can work out the generators and relations for the corresponding big algebras explicitly. In particular, one can also compute explicitly  $\calB^{2\omega_2}(\so_5)\subset \calC^{2\omega_2}(\so_5)$ the big algebra of the adjoint representation of $\SO_5$. We can obtain this adjoint representation by restricting the representation $\rho_{\omega_2}$ of $\SL_5$ to the subgroup $\SO_5\subset \SL_5$. This way we also have a commutative subalgebra $\calB^{\omega_2}(\sl_5)\otimes_{H^{2*}_{\SL_5}} H^{2*}_{\SO_5}\subset \calC^{2\omega_2}(\so_5).$ Both subalgebras of  $\calC^{2\omega_2}(\so_5)$ satisfy  properties 2.,3. and 4. in Theorem~\ref{big} but can be shown to be non-isomorphic. This shows that the big algebra $\calB^{2\omega_2}(\so_5)\subset \calC^{2\omega_2}(\so_5)$ is not  uniquely determined by these properties.
\end{remark}

\subsection{Twining big algebras}
\label{twiningsection}

For a connected semisimple complex Lie group $\G$ let $\sigma:\G\to \G$ be a distinguished automorphism, i.e. one which fixes a pinning. In particular, it is induced from an automorphism, also denoted $\sigma$, of the Dynkin diagram. Examples for the symmetric pair $(\G,\G^\sigma)$ are $(\SL_{2n+1},\SO_{2n+1}),$ $(\SL_{2n},\Sp_n),$ $(\SO_{2n},\SO_{2n-1}),$  $(\rm{PSO}_8,\rm{G}_2)$ or $(\rm{E}_6,\rm{F}_4).$ Except for the order three $\sigma$ in the case $(\rm{PSO}_8,\rm{G}_2)$ the automorphism $\sigma$ is an involution. 

The Dynkin diagram automorphism $\sigma$ induces a distinguished automorphism $\sigma: \G^\vee\to \G^\vee$ of the Langlands dual. Define the endoscopy group $\G_\sigma=((\G^\vee)^\sigma_0)^\vee $. Such a $\sigma$ will induce an automorphism of the Feigin-Frenkel center, the Gaudin algebra and the universal big algebra, and in turn for $\mu\in \Lambda^+(\G)^\sigma$ on the big algebra $\sigma:\calB^\mu\to \calB^\mu$. Decompose $\calB^\mu=\oplus_{\kappa\in \widehat{\langle \sigma\rangle}}(\calB^\mu)_\kappa$ according to characters of the cyclic group $\langle \sigma \rangle\subset Aut(\G)$.  Define the coinvariant algebra $\calB^\mu_\sigma:=\calB^\mu/(\oplus_{1\neq \kappa \in \widehat{\langle\sigma \rangle} } (\calB^\mu)_\kappa)$, which computes the ring of functions of the fixed point scheme: $\calB^\mu_\sigma\cong \C[\Spec(\calB^\mu)^\sigma]$. We have the following\footnote{A proof of this conjecture appeared in \cite{zveryk}.}  
\begin{conjecture}\label{twiningconj} For $\mu\in \Lambda^+(\G_\sigma)$  also denote the corresponding dominant weight by $\mu\in\Lambda^+(\G)^\sigma$. Then \begin{align}\calB^\mu_\sigma (\g)\cong \calB^\mu(\g_\sigma)\label{twining}.\end{align}
    \end{conjecture}

    The main motivation for the conjecture was that it is compatible with Jantzen's twining character formula. Namely take $\lambda\in \Lambda^+(\G_\sigma)$ and the corresponding $\lambda\in \Lambda^+(\G)^\sigma$. The weight space $V^\mu_\lambda(\G)$ of the $\G$-representation will inherit an action $\sigma:V^\mu_\lambda(\G)\to V^\mu_\lambda(\G)$, which combined with the induced action in the big algebra $\sigma:\calB^\mu\to \calB^\mu$ will yield an automorphism of the multiplicity algebra $Q^\mu_\lambda(\g)$. Then  we expect \eqref{twining} implies that $Q^\mu_\lambda(\g)_\sigma=Q^\mu_\lambda(\g_\sigma)$ and  $\dim(Q^\mu_\lambda(\g)_\sigma)=\tr(\sigma:Q^\mu_\lambda(\g)\to Q^\mu_\lambda(\g))$, when the trace is non-zero. In this case we get that 
        $\tr(\sigma:V^\mu_\lambda(\G)\to V^\mu_\lambda(\G))=\tr(\sigma:Q^\mu_\lambda(\g)\to Q^\mu_\lambda(\g))=\dim(Q^\mu_\lambda(\g_\sigma))=\dim(V^\mu_\lambda(\G_\sigma))$, which is Jantzen's twining formula \cite[Satz 9]{jantzen}. 

 Geometrically the result should  follow from the induced action $\sigma:\Gr^\mu(\G^\vee)\to \Gr^\mu(\G^\vee)$ for $\mu\in \Lambda^+(\G)^\sigma$. In fact, the first check on the conjecture is when $V^\mu(\G)$ is a $\sigma$-invariant minuscule representation. When $\mu=\omega_n\in \Lambda^+(\SL_{2n})$ then $\sigma(\mu)=\mu$ and the corresponding cominuscule flag variety $\Gr^{n\omega_1}(\PGL_{2n})\cong \Gr(n,\C^{2n})$ is the Grassmannian of $n$-planes in $\C^{2n}$. The action of $\sigma$ on $\Gr(n,\C^{2n})$ is given by $\sigma(V):=\rm{ann}(\omega(V))$, where $\omega:\C^{2n}\to (\C^{2n})^*$ is a symplectic form.  Thus we see that $\Gr(n,\C^{2n})^\sigma\cong \rm{L}\Gr(n,\C^{2n})\cong \Gr^{\omega_n}(\PSp_{2n})$ is the Lagrangian Grassmannian. As $\Gr(n,\C^{2n})$ is $\PGL_{2n}$-regular and $\rm{L}\Gr(n,\C^{2n})$ is $\PSp_{2n}$-regular, from \cite[Theorem 1.3]{hausel-rychlewicz} we can deduce that $$\calB^{\omega_n}(\sl_{2n})_\sigma\cong H^{2*}_{\PGL_{2n}}(\Gr(n,\C^{2n}))_\sigma\cong H^{2*}_{\PGL_{2n}^\sigma}(\Gr(n,\C^{2n})^\sigma)\\ \cong H^{2*}_{\PSp_{2n}}(\rm{L}\Gr(n,\C^{2n}))\cong \calB^{\omega_n}(\so_{2n+1}).$$

 Finally we note that in the example $\G= \PGL_3$ and $\G_\sigma= \SL_2$  the weight $\omega_1\in \Lambda^+(\SL_2)$ corresponds to $\omega_1+\omega_2\in \Lambda^+(\PGL_3)^\sigma$. Then we have $\sigma(M_1)=M_1$, $\sigma(N_1)=-N_1$ and $\sigma(M_2)=-M_2$ and so the corresponding $\calB^{\omega_1+\omega_2}(\sl_3)_\sigma\cong \calB^{\omega_1}(\sl_2)$ can be seen in the first picture of Figure~\ref{octetfigure}. Namely, the fixed point scheme of  $\sigma$ on $\Spec(\calB^{\omega_1+\omega_2}(\sl_3))$ is the common parabola of the big skeleton shared with the medium skeleton, where $N_1=M_2=0$. 

\subsection{Mirror symmetry and big spectral curves}

Big algebras first appeared in \cite{hausel-madrid} in connection with mirror symmetry \cite{hausel-hitchin,hausel-ICM}. 
They were needed to endow the universal $\G$-Higgs bundle in an irreducible representation with the structure of a bundle of algebras  along the Hitchin section. 
Turning the logic back, one can use the big algebras $\calB^\mu$ to define a bundle of algebras  on the $\G$-Higgs bundle  in the irreducible representation $V^\mu$ along the Hitchin section, yielding {\em big spectral curves} $C^\mu\subset \oplus_{k=1}^{rank(\G)}\oplus_{0<i<k}K^{d_k-i}$ living in the total space of direct sum of line bundles $K^i$ for each degree $i$ generator of the big algebra. In turn, for any $\G$-Higgs bundle one can construct a big algebra of big Higgs fields in any irreducible representation $V^\mu$, which will yield a rank $1$ sheaf on the corresponding big spectral curve $C^\mu$. We expect a full theory of BNR correspondences for each big spectral curve, bridging the usual spectral curves in \cite{hitchin}
with the cameral covers in \cite{donagi}.

Finally, we expect that the geometric description of the quantum big algebras $\calG^\mu$ in \cite{feigin-frenkel-toledano} as rings of functions on certain spaces of opers, and the description \cite{hausel-madrid} of the big algebras $\calB^\mu$ as rings of functions on upward flows in the Hitchin system could be unified as a description of the $\hbar$-quantum big algebras $\calG^\mu_\hbar$ on  upward flows in $\M_{Hodge}$, the moduli space of $\hbar$-connections. 

Details of the proofs of the results in this paper, and detailed study of the examples mentioned above will appear elsewhere. 
\vskip.2cm
{\noindent \bf Acknowledgements.} We thank Nigel Hitchin for discussions and 
 the joint projects \cite{hausel-hitchin,hausel-ICM} this paper has grown out from. We thank Vladyslav Zveryk for collaboration 
 on Theorem~\ref{explicit} and on the corresponding Magma code which implements big algebras.  We thank Hiraku Nakajima for discussions and pointing out Theorem~\ref{geometric}.\ref{nakajima}, a result generalising our original observation in the $\hbar=0$ case. Special thanks go to Leonid Rybnikov for patiently explaining his works, in particular \cite{feigin-frenkel-rybnikov} crucial to Theorem~\ref{big}. We thank Michel Brion,
Michael Finkelberg,  Oscar Garc\'ia-Prada, Jakub L\"owit, Joel Kamnitzer, Friedrich Knop, Michael McBreen, Anton Mellit, Takuro Mochizuki, Shon Ng\^o,  Kamil Rychlewicz, Shiyu Shen, Leslie Spencer, Bal\'azs Szendr\Horig{o}i, Andr\'as Szenes and Oksana Yakimova for comments and discussions. Kamil Rychlewicz and Daniel Bedats helped with the Mathematica files for the figures, and we used the SM\_isospin Tikz package of Izaak Neutelings for drawing the baryon multiplets. We thank the referees for many useful comments. We acknowledge funding from FWF grant "Geometry of the tip of the global nilpotent cone" no.number P 35847.

\bibliographystyle{acm}

\end{document}